\documentclass[12pt]{amsart}  
\usepackage{amsmath,amsthm,amsfonts,amssymb} 
\numberwithin{equation}{section} 
\newtheorem{thm}[equation]{Theorem}\newtheorem{cor}[equation]{Corollary} 
\newtheorem{lem}[equation]{Lemma}  
\newtheorem{example}[equation]{Example}
\newtheorem{remark}[equation]{Remark}
\theoremstyle{definition}

\DeclareMathOperator{\codim}{codim}

\DeclareMathOperator{\diag}{diag}

\DeclareMathOperator{\Ima}{im}
\DeclareMathOperator{\Ker}{ker}

\newcommand{\C}{\mathbb{C}}
\newcommand{\coh}{{\rm H}}
\newcommand{\DOT}{\setlength{\unitlength}{1pt}\begin{picture}(2.5,2)
                  (1,1)\put(2,3.5){\circle*{2}}\end{picture}}
\newcommand{\HH}{{\rm HH}}

\newcommand{\ot}{\otimes}

\newcommand{\Wedge}{\textstyle\bigwedge}
\newcommand{\Z}{\mathbb {Z}}

\begin{document}

\title[TWISTED GRADED HECKE ALGEBRAS]
{TWISTED GRADED HECKE ALGEBRAS}

\author{Sarah Witherspoon}
\address{Department of Mathematics\\
 Texas A\&M University \\ College Station, Texas 77843, USA} 
\email{sjw@math.tamu.edu} 
\thanks{Partially supported by NSF grants \#{}DMS--0422506, 0443476.}

\date{June 7, 2005}


\begin{abstract}

We generalize graded Hecke algebras to include a twisting two-cocycle
for the associated finite group. 
We give examples where the parameter spaces of the resulting 
{\em twisted} graded Hecke algebras are larger
than that of the graded Hecke algebras. We prove that twisted
graded Hecke algebras are particular types of deformations of a
crossed product.

\end{abstract}

\maketitle  

\section{Introduction}

Drinfel'd defined graded Hecke algebras for any finite subgroup $G$ of
$GL(V)$ \cite{drinfeld86}. This definition was shown by Ram and
Shepler \cite{ram-shepler03} to generalize the definition of graded versions
of affine Hecke algebras for real reflection groups given by Lusztig,
who was motivated by questions in representation theory
\cite{lusztig88,lusztig89}. In the same paper Ram and Shepler
classified all graded Hecke algebras for complex reflection 
groups. They found many nontrivial graded Hecke algebras, but also 
showed that some groups (such as $G(r,1,n)=\Z/r\Z\wr S_n$ when
$r>2, \ n>3$) have no nontrivial graded Hecke algebras.
Etingof and Ginzburg showed that in the case of a finite symplectic
group, the symplectic form itself arises naturally in the structure
of any associated graded Hecke algebra, always yielding nontrivial
examples, which they called symplectic reflection algebras
\cite{etingof-ginzburg02}. Many authors have
studied representations of graded Hecke algebras, their subalgebras 
generated by certain idempotents, and connections to the
geometry of the corresponding orbifolds $V/G$.
(See for example \cite{gordon-smith04,kriloff-ram02,lusztig95,lusztig02}.)

In this note we generalize these constructions by incorporating a
two-cocycle $\alpha$ that represents an element of the cohomology group
$\coh^2(G,\C^{\times})$. Chmutova introduced such a two-cocycle
for symplectic groups, showing that it appears naturally in
symplectic reflection algebras arising from nonfaithful
group representations to {\em Sp}$(V)$ \cite{chmutova}.
Here our motivation comes from numerous papers on orbifolds in which
such a cocycle, called {\em discrete torsion}, appears 
(\cite{adem-ruan03,caldararu-giaquinto-witherspoon04,vafa-witten95}
are just a few), as well as our finding that for some groups
there is a nontrivial {\em twisted} graded Hecke algebra even if there is
no nontrivial graded Hecke algebra (such as Example \ref{ex:1} below).
We adapt the direct linear-algebraic approach of Ram and Shepler to derive
criteria for existence of such nontrivial twisted graded Hecke algebras
(Theorem \ref{thm:necsuff}, Corollary \ref{cor:parameters}). 
We give examples for which the parameter spaces for twisted graded
Hecke algebras are larger than that for graded Hecke algebras:
$S_n\subset GL(\C^n)$ via the permutation representation
(Example \ref{ex:3}), as well as
$(\Z/\ell Z)^{n-1}\subset SL(\C^n)$ (Example \ref{ex:1}).
On the other hand, if $G$ is symplectic, this is not the case
(Example \ref{ex:2}).
It would be interesting to do a more thorough analysis of the
parameter spaces of twisted graded Hecke algebras for various types
of groups.

In Section 3 we
show that twisted graded Hecke algebras are precisely particular types of
deformations of crossed product algebras 
$S(V)\#_{\alpha}G$ (Theorem \ref{thm:equivalence}).
We use this
idea to show that Ram and Shepler's isomorphism between the different
definitions of (untwisted) graded Hecke algebra given by Drinfel'd
and Lusztig arises as an equivalence of deformations whose infinitesimals
are thus cohomologous. 
This also puts previous results into a larger context: There are in
general many more infinitesimal deformations (that is, Hochschild
two-cocycles) of $S(V)\#_{\alpha}G$ than those that lift to
deformations that are (twisted) graded Hecke algebras.
In some cases these other infinitesimals also lift to deformations
of $S(V)\#_{\alpha}G$: See Remark \ref{deg1} herein and
\cite{caldararu-giaquinto-witherspoon04,witherspoon}.
We are not aware of a general result
regarding such deformations.

We thank A.\ Ram for several stimulating conversations. It was a
question of his that ultimately led to this work as well as
to insight into other related projects.
We thank R.-O.\ Buchweitz for explaining a computation of the
relevant Hochschild cohomology to us.

Throughout, we will work over the complex numbers $\C$, so that
$\ot = \ot_{\C}$ unless otherwise indicated.

\section{Twisted graded Hecke algebras}

Our approach to graded Hecke algebras will be through 
crossed products (generalizations of skew group algebras), 
so we begin by summarizing this construction. 
Our definition below of (twisted) graded
Hecke algebras involves a parameter $t$, which may be taken to
be an indeterminate or any complex number.
In Section 3 we will assume $t$ is an indeterminate,
as we will discuss deformations over the polynomial ring $\C[t]$.
Specializing $t$ to any nonzero complex number
results in a definition of graded Hecke algebra
equivalent to those in the literature.

Let $V=\C^n$, $G$ a finite subgroup of $GL(V)$, and $\alpha:
G\times G\rightarrow \C^{\times}$ a two-cocycle, that is
\begin{equation}\label{eqn:alpha}
   \alpha(g,h)\alpha(gh,k)=\alpha(h,k)\alpha(g,hk)
\end{equation}
for all $g,h,k\in G$. The action of $G$ on $V$ induces actions
by algebra automorphisms on the tensor algebra $T(V)$ and the
symmetric algebra $S(V)$. 

Suppose $S$ is any associative $\C$-algebra with an action of $G$
by automorphisms, for example $S=T(V)$ or $S=S(V)$. Then we may
form the {\em crossed product algebra} $S\#_{\alpha}G$, which is
$S\otimes \C G$ as a vector space, and has multiplication
$$
  (r\ot g) (s\ot h) = \alpha(g,h) r(g\cdot s)\ot gh
$$
for all $r,s\in S$ and $g,h\in G$. This makes $S\#_{\alpha}G$
an associative algebra as $\alpha$ is a two-cocycle. Note that
$S$ is a subalgebra of $S\#_{\alpha}G$, via $\displaystyle{
S\stackrel{\sim}{\longrightarrow} S\ot 1}$, and that the 
subalgebra $1\ot \C G$ of $S\#_{\alpha} G$ is also known as the
{\em twisted group algebra} $\C^{\alpha}G$. We accordingly 
abbreviate $s\ot g$ by $s\overline{g}$ ($s\in S, \ g\in G$).
The action of $G$ on $S$ becomes an inner action on $S\#_{\alpha}G$,
namely $s\mapsto \overline{g} s (\overline{g})^{-1}$ for all $s\in S,
\ g\in G$.

The two-cocycle $\alpha$ is a {\em coboundary} if there is some
function $\beta: G\rightarrow \C^{\times}$ such that
\begin{equation}\label{eqn:coboundary}
  \alpha(g,h) = \beta(g)\beta(h)\beta(gh)^{-1}
\end{equation}
for all $g,h\in G$. The set of two-cocycles modulo coboundaries
forms an abelian group under pointwise multiplication, called the
{\em Schur multiplier} of $G$, and denoted $\coh ^2(G,\C^{\times})$.
If two cocycles differ by a coboundary, that is if they are
{\em cohomologous}, the corresponding crossed products are isomorphic.

Replacing $\alpha$ by a cohomologous cocycle if necessary,
we may assume that $\alpha$ is normalized so that $\alpha(1,g) =
\alpha(g,1) = 1$ for all $g\in G$. It follows 
that $\overline{1}$ is the multiplicative identity for $\C^{\alpha}G$.
Thus $(\overline{g})^{-1}=\alpha^{-1}(g,g^{-1})\overline{g^{-1}}=
\alpha^{-1}(g^{-1},g)\overline{g^{-1}}$ for all $g\in G$.

For each $g\in G$, choose a skew-symmetric bilinear form
$a_g:V\times V\rightarrow \C$ (arbitrary for now, and possibly 0).
Let $t$ be an indeterminate and
extend scalars, $S\#_{\alpha}G[t]:= \C[t]\ot (S\#_{\alpha}G)$.
Let $A$ be the quotient of $T(V)\#_{\alpha}G[t]$ by the 
(two-sided) ideal $I[t]$ generated by all 
\begin{equation}\label{eqn:bracket}
  [v,w] - \sum_{g\in G} a_g(v,w) t \overline{g}
\end{equation}
where $v,w\in V$ and $[v,w]=vw-wv$. (We have omitted tensor symbols
in elements of $T(V)$ for brevity.) 
Note that $A$ is {\em additively}
isomorphic to a quotient of $S(V)\#_{\alpha}G [t]$ 
via any choice of linear section $S(V)\hookrightarrow T(V)$
of the canonical projection of $T(V)$ onto $S(V)$. 
This quotient is proper if and only if there
is more than one way to rearrange factors in a word by applying the
relations (\ref{eqn:bracket}). We say that $A$ is a {\em twisted
graded Hecke algebra} if the above map yields $A\cong S(V)\ot \C G[t]$ 
as vector spaces over $\C[t]$. 
Equivalently, if we assign degree 1 to elements of $V$ and degree 0
to elements of $G$, then $A$ is a filtered algebra over $\C [t]$
and is a twisted graded Hecke algebra in case the associated
graded algebra ${\rm gr}A$ is isomorphic to $S(V)\#_{\alpha}\C G[t]$.

\begin{remark}{\em
Let $\alpha$ be the trivial two-cocycle $\alpha(g,h)=1$ for
all $g,h\in G$, and let $t=1$. Then $A$
becomes the {\em graded Hecke algebra} of \cite{drinfeld86,ram-shepler03}.
}
\end{remark}

We will use the techniques of Ram and Shepler \cite{ram-shepler03}
to determine the conditions on the bilinear forms $a_g$ under
which $A$ is a twisted graded Hecke algebra. Let $h\in G$, and 
use (\ref{eqn:bracket}) to obtain two expressions involving
$[h^{-1}\cdot v,h^{-1}\cdot w]$. By direct substitution:
$[h^{-1}\cdot v, h^{-1}\cdot w]=\sum_{g\in G} a_g(h^{-1}\cdot v,
h^{-1}\cdot w) t\overline{g}$. By conjugating (\ref{eqn:bracket}) by 
$(\overline{h})^{-1}$:
\begin{eqnarray*}
  [h^{-1}\cdot v,h^{-1}\cdot w] &=& \sum_{g\in G} a_g(v,w) t
  (\overline{h})^{-1} \overline{g}\overline{h}\\
  &=& \sum_{g\in G}a_g(v,w)\alpha^{-1}(h,h^{-1})\alpha(g,h)
  \alpha(h^{-1},gh) t\overline{h^{-1}gh}.
\end{eqnarray*} 
Both must be in the ideal $I[t]$, and the
assumed additive isomorphism
$A\cong S(V)\ot \C G [t]$ allows us to equate coefficients of
each $\overline{g}$ ($g\in G$) in these two expressions.
Thus we have $a_{h^{-1}gh}(h^{-1}\cdot v,h^{-1}\cdot w) =
\alpha^{-1}(h,h^{-1})\alpha(g,h)\alpha(h^{-1},gh)a_g(v,w)$, and
replacing $v,w$ by $h\cdot v,h\cdot w$, we
obtain the first equation below.
Again the isomorphism $A\cong S(V)\ot \C G[t]$ implies that in any
expression $wvu\in V^{\ot 3}$, three applications of (\ref{eqn:bracket})
in any order to obtain $uvw$ plus an element in $tA$ must yield the
same element of $A$. Comparison results in the second equation below. 
Thus necessary conditions for $A$ to be a twisted graded Hecke algebra
are:
\begin{equation}\label{eqn:conjugate}
  a_{h^{-1}gh}(v,w)=\alpha^{-1}(h,h^{-1})
  \alpha(g,h)\alpha(h^{-1},gh)a_g(h\cdot v,h\cdot w)
\end{equation}
\begin{equation}\label{eqn:jacobi}
  a_g(u,v) (g\cdot w - w) + a_g(v,w)(g\cdot u-u)
      +a_g(w,u) (g\cdot v-v) = 0
\end{equation}
for all $g,h\in G$ and $u,v,w\in V$. These equations will be used in
the proofs of the following lemma and theorem. 

Not only are (\ref{eqn:conjugate}) and (\ref{eqn:jacobi}) necessary
for $A$ to be a twisted graded Hecke algebra, but
they are also sufficient: The relations (\ref{eqn:bracket})
allow for rearrangement of any expression in $A$ to a particular
form identified with an element of $S(V)\ot \C G[t]$, and the 
relations (\ref{eqn:conjugate})
and (\ref{eqn:jacobi}) imply that such a form is unique.
That is, (\ref{eqn:conjugate}) is equivalent to uniqueness of the
canonical form of $\overline{h}vu$, and (\ref{eqn:jacobi})
is equivalent to uniqueness of the canonical form of
$wvu$ ($u,v,w\in V$, $h\in G$), as demonstrated in 
\cite{ram-shepler03} for the case
$\alpha =1$ (in the text above Lemma 1.5).
The uniqueness of the form of $wvu$ is equivalent to the Jacobi identity 
in $A$,
\begin{equation}\label{eqn:jacobi2}
[u, [v,w]]+[v,[w,u]]+[w,[u,v]]=0,
\end{equation}
which is another way to express the condition (\ref{eqn:jacobi}).
For the uniqueness of the canonical form of a monomial having
more than three factors, it is helpful to note that replacing any
of $u,v,w$ in the left side of (\ref{eqn:jacobi2}) by an element
of degree larger than 1 results in an element in the ideal generated
by the left side of (\ref{eqn:jacobi2}). (For more details, see
\cite{ram-shepler03}.)

Let $( \ , \ )$ be any $G$-invariant nondegenerate
Hermitian form on $V$. Define orthogonal complements of subspaces
of $V$ via this form.
Let $g\in G$ and $V^g=\{v\in V\mid g\cdot v=v\}$.
We will need the observation that
$(V^g)^{\perp} =\Ima (g - 1)$, which follows from the
standard facts $V^g=\ker (g-1)$, $\Ima (g-1)\subset (V^g)^{\perp}$,
and $\dim (\Ima (g-1))=\dim (\ker (g-1)^{\perp})$.
For any $h\in G$, denote by
$h^{\perp} : (V^g)^{\perp}\rightarrow (V^g)^{\perp}$ the composition
of the linear maps $h: (V^g)^{\perp} \hookrightarrow V$,
$V\rightarrow V/V^g$, and a choice of isomorphism
$\displaystyle{V/V^g \stackrel{\sim}{\rightarrow} (V^g)^{\perp}}$.

\begin{lem}\label{lem:necessary}
Suppose that $A$ is a twisted graded Hecke algebra for $G$ defined
by skew-symmetric bilinear forms $a_g$ ($g\in G$) on $V$. For each
$g\neq 1$, if $a_g\neq 0$ then $\Ker a_g = V^g$, $\codim(V^g)=2$, and
$$
  a_g(h\cdot v, h\cdot w)=\det(h^{\perp}) a_g(v,w)
$$
for all $h\in G$ and $v,w\in (V^g)^{\perp}$.
\end{lem}

This proof is essentially the same as in Ram and Shepler 
\cite[Lemma 1.8]{ram-shepler03},
with appropriate adjustments made for the possibly
nontrivial two-cocycle $\alpha$.

\begin{proof}
Let $g\neq 1$ with $a_g\neq 0$. Let $v\in V$ and $w\in V^g$.
If $v$ is also in $V^g$, then $a_g(v,w)(g\cdot u-u)=0$ for all $u\in V$ by
(\ref{eqn:jacobi}). As $g\neq 1$, there is some $u\in V$ with
$g\cdot u\neq u$, so $a_g(v,w)=0$. If $v\not\in V^g$, let
$u=\sum_{k=1}^r g^k\cdot v$, where $r$ is the order of $g$, so
that $u\in V^g$ (possibly $u=0$). 
As before, this implies $a_g(u,w)=0$, and we
may rewrite this equation using the definition of $u$ and
(\ref{eqn:conjugate}) as
\begin{eqnarray*}
  0=a_g(u,w) &=& \sum_{k=1}^r a_g(g^k\cdot v,w)\\
             &=& \sum_{k=1}^r \alpha^{-1}(g^{-k},g^k)\alpha(g,g^{-k})
              \alpha(g^k,g^{1-k})a_g(v,w).
\end{eqnarray*}
The restriction of $\alpha$ to the cyclic subgroup $\langle g\rangle$
generated by $g$ is a coboundary since the Schur multiplier of a cyclic
group is trivial \cite[Prop.\ 1.1]{karpilovsky85}. Therefore there is
a function $\beta:\langle g\rangle\times\langle g\rangle\rightarrow
\C^{\times}$ such that $\alpha(h,l)=\beta(h)\beta(l)\beta^{-1}(hl)$
for all $h,l\in\langle g\rangle$. It follows that the above coefficients
of $a_g(v,w)$ are
\begin{eqnarray*}
 &&\hspace{-.5in} \alpha^{-1}(g^{-k},g^k)\alpha(g,g^{-k})\alpha(g^k,g^{1-k})\\
   \hspace{-3in}&=&\beta^{-1}(g^{-k})\beta^{-1}(g^k)
\beta(1)\beta(g)\beta(g^{-k})
  \beta^{-1}(g^{1-k})\beta(g^k)\beta(g^{1-k})\beta^{-1}(g)\\
  \hspace{-3in}&=& \beta(1),
\end{eqnarray*}
that is they are independent of $k$. The above calculation reduces to
$0=a_g(u,w)=r\beta(1)a_g(v,w)$. As $r\beta(1)\neq 0$, this forces
$a_g(v,w)=0$.
Consequently $V^g\subset \Ker a_g$, and so $(\Ker a_g)^{\perp} 
\subset (V^g)^{\perp}$.

As $a_g$ is nonzero and skew-symmetric, we have 
$\dim((\Ker a_g)^{\perp})=\codim (\Ker a_g)\geq 2$. Let $u,v$ be two linearly
independent elements of $(\Ker a_g)^{\perp}\subset (V^g)^{\perp}$
with $a_g(u,v)\neq 0$, so that in particular $g\cdot u - u\neq 0$
and $g\cdot v - v\neq 0$. Let $w$ be any element of $(V^g)^{\perp}
=\Ima (g-1)$, and write $w=g\cdot w'-w'$ for some $w'\in V$.
By (\ref{eqn:jacobi}), $w$ is a linear combination of $g\cdot u-u$
and $g\cdot v -v$, two elements of $(V^g)^{\perp}$. This implies
$\dim(\Ker a_g)^{\perp}=\dim(V^g)^{\perp} =2$, so that 
$V^g=\Ker a_g$ and $\codim(V^g)=2$.

Finally, let $v,w\in (V^g)^{\perp}$, $h\in G$, and
\begin{eqnarray*}
  h\cdot v &=& a_{11} v + a_{21}w +x\\
  h\cdot w &=& a_{12}v+a_{22}w+x'
\end{eqnarray*}
where $x,x'\in V^g$ and $a_{ij}$ are scalars. Applying 
(\ref{eqn:conjugate}), we have $a_{h^{-1}gh}(v,w)=\alpha^{-1}(h,h^{-1})
\alpha(g,h)\alpha(h^{-1},gh)a_g(h\cdot v,h\cdot w)$ on the one hand, while
evaluating via the above equations yields
\begin{equation}\begin{array}{rcl}
  a_{h^{-1}gh}(v,w) &=& \alpha^{-1}(h,h^{-1})
                  \alpha(g,h)\alpha(h^{-1},gh)a_g(a_{11}v
                   +a_{21}w+x, a_{12}v+a_{22}w +x')\\
   &=& \alpha^{-1}(h,h^{-1})\alpha(g,h)\alpha(h^{-1},gh)
     (a_{11}a_{22}-a_{12}a_{21}) a_g(v,w)\\
   &=& \alpha^{-1}(h,h^{-1})\alpha(g,h) 
     \alpha(h^{-1},gh)\det(h^{\perp})a_g(v,w).
\end{array}\label{eqn:ahgh}
\end{equation}
Equating the two expressions for $a_{h^{-1}gh}(v,w)$ gives
$a_g(h\cdot v,h\cdot w)=\det(h^{\perp}) a_g(v,w)$, as desired.
\end{proof}

Now suppose $g\neq 1$, $a_g\neq 0$, and $h\in C(g)=\{h\in G\mid hg=gh\}$.
Then (\ref{eqn:ahgh}) is equivalent to
$\det(h^{\perp})=\alpha(h,h^{-1})\alpha^{-1}(g,h)\alpha^{-1}
(h^{-1},gh)$. Applying (\ref{eqn:alpha}) to the triple $h,h^{-1},gh$,
we find that $\alpha(h,h^{-1})\alpha^{-1}(h^{-1},gh)=\alpha(h,h^{-1}gh)
=\alpha(h,g)$ when $h\in C(g)$,
so that $\det(h^{\perp})=\alpha(h,g)\alpha^{-1}(g,h)$ for all $h\in C(g)$.
This condition is independent of the choice of $\alpha$ in a given
coset modulo coboundaries since coboundaries are symmetric on commuting pairs
as is evident from (\ref{eqn:coboundary}). This is as expected since the
determinant function is independent of such choices. Note that
as $g\in C(g)$, (\ref{eqn:ahgh}) implies $\det(g^{\perp})=1$, that is
$g\in SL(V)$. Further, in case $g$ is $\alpha$-{\em regular}, 
that is $\alpha(g,h)=\alpha(h,g)$ for all $h\in C(g)$, this determinant
condition is simply $\det(h^{\perp})=1$ for all $h\in C(g)$. For
nonregular elements $g$, this condition is different from that in the case of
the trivial cocycle, leading to new examples such as Examples \ref{ex:1}
and \ref{ex:3} below.

\begin{thm}\label{thm:necsuff}
Let $G$ be a finite subgroup of $GL(V)$, $\alpha:G\times G
\rightarrow \C^{\times}$ a normalized two-cocycle, and 
$g\in G-\{1\}$. There is a twisted graded Hecke algebra $A$ with
$a_g\neq 0$ if, and only if, $\Ker a_g=V^g$, $\codim (V^g)=2$, and
\begin{equation}\label{eqn:condition}
\det(h^{\perp})=\alpha(h,g)\alpha^{-1}(g,h)
\end{equation}
 for all $h\in C(g)$.
\end{thm}

Again the proof is similar to that of Ram and Shepler 
\cite[Thm.\ 1.9]{ram-shepler03}
in the untwisted case. See also the paper of Etingof and Ginzburg
\cite{etingof-ginzburg02}, who used a criterion of Braverman and
Gaitsgory adapted to Koszul algebras over $\C G$
\cite{beilinson-ginzburg-soergel96,braverman-gaitsgory96}.

\begin{proof}
If $A$ is a twisted graded Hecke algebra, Lemma \ref{lem:necessary}
and subsequent comments show the given conditions hold.

Conversely, suppose the stated conditions hold for $g$. Up to a
scalar multiple, there is a unique skew-symmetric form on $V$ that is
nondegenerate on $(V^g)^{\perp}$ and has kernel $V^g$. Fix such a
form $a_g$, and let
\begin{equation}\label{eqn:ak}
  a_k(v,w)=\left\{\begin{array}{cc}
    \alpha^{-1}(h,h^{-1})\alpha(g,h)\alpha(h^{-1},gh) a_g(h\cdot v,h\cdot w), &
    \mbox{ if } k=h^{-1}gh\\
    0, & \mbox{ otherwise}\end{array}\right.
\end{equation}
for all $v,w\in V$. In order for (\ref{eqn:conjugate}) to hold for all
pairs of group elements, we must check that the definition of $a_k$
is independent of the choice of representative from the conjugacy class
of $g$. Suppose that $k=l^{-1}h^{-1}ghl$, so that there are two ways 
to define $a_k$. One way yields
$$
  a_k(v,w)=\alpha^{-1}(hl,l^{-1}h^{-1})\alpha(g,hl)
  \alpha(l^{-1}h^{-1},ghl)a_g(hl\cdot v,hl\cdot w),
$$
and the other yields
$$
   a_k(v,w)\hspace{5.1in}
$$
$$
=\alpha^{-1}(l,l^{-1})\alpha(h^{-1}gh,l)\alpha(l^{-1},h^{-1}ghl)
   \alpha^{-1}(h,h^{-1})\alpha(g,h)
   \alpha(h^{-1},gh)a_g(hl\cdot v,hl\cdot w).
$$
These two expressions for $a_k(v,w)$ are indeed equal, as follows from
five applications of (\ref{eqn:alpha}), to the triple $h^{-1},gh,l$, 
to $l^{-1},h^{-1},ghl$, to $g,h,l$, to $hl,l^{-1},h^{-1}$, and
to $h,l,l^{-1}$.
By the discussion before Theorem \ref{thm:necsuff}, if $h^{-1}gh=g$, then
(\ref{eqn:ak}) coincides with $a_g(v,w)$.
Therefore the definition of $a_k$ is independent of the choice of 
representative from the conjugacy class of $g$, and
(\ref{eqn:conjugate}) holds. The
argument of \cite[Lemma 1.8(b)]{ram-shepler03} applies without change
to show that (\ref{eqn:jacobi}) holds as well: Their key observation
is that (\ref{eqn:jacobi}) holds trivially if any one of $u,v,w$
is in $V^g$, while $\dim(V^g)^{\perp}=2$ implies that any three
elements $u,v,w$ of $(V^g)^{\perp}$ must be linearly dependent.
Substituting a linear dependence relation into the left side of 
(\ref{eqn:jacobi}) yields 0 after some manipulation. Therefore $A$,
defined by the $\{a_k\}$ in (\ref{eqn:ak}), is a twisted graded
Hecke algebra.
\end{proof}

In summary, under the conditions in the theorem, $a_g$ is 
determined by its value $a_g(v,w)$ on a basis $v,w$ of $(V^g)^{\perp}$,
and for each $h\in G$, $a_{h^{-1}gh}(v,w)$ is given by (\ref{eqn:ak}).

\begin{cor}\label{cor:parameters}
Let $d$ be the number of conjugacy classes of $g\in G$ such that
$\codim V^g =2$ and $\det(h^{\perp})=\alpha(h,g)\alpha^{-1}(g,h)$
for all $h\in C(g)$. The sets $\{a_g\}_{g\in G}$ corresponding to
twisted graded Hecke algebras $A$ form a vector space of dimension
$d+\dim (\Wedge^2V)^G$.
\end{cor}

\begin{proof}
By (\ref{eqn:conjugate}) and (\ref{eqn:jacobi}), the only
condition on $a_1$ is $a_1(v,w)=a_1(h\cdot v,h\cdot w)$ for all
$h\in G$ since $\alpha(h,h^{-1})=\alpha(h^{-1},h)$, 
that is $a_1$ is a $G$-invariant element of $(\Wedge^2V)^*$.
For each conjugacy class of elements $g\in G$ with $\codim(V^g)=2$
and $\det(h^{\perp})=\alpha(h,g)\alpha^{-1}(g,h)$ for all $h\in C(g)$,
there are skew-symmetric forms $a_g$ with kernel $V^g$, 
determined by a single form $a_g$ (unique up to scalar multiple)
for a representative $g$, and given by (\ref{eqn:ak}).
\end{proof}

\begin{remark}\label{deg1}{\em
The results of this section apply if we relax one of the conditions on
the forms $a_g$ a little: Allow $a_g$ to take values in $\C\oplus V^g$, 
that is $g$-invariant polynomials of degree at most 1. Again let $A$
be the quotient of $T(V)\#_{\alpha} G[t]$ by the ideal $I[t]$
generated by all expressions of the form (\ref{eqn:bracket}). As the
image of each form $a_g$ has degree at most 1, the techniques of Ram and
Shepler \cite{ram-shepler03} still apply to yield
$$
  a_{h^{-1}gh}(v,w)=\alpha^{-1}(h,h^{-1})\alpha(g,h)\alpha(h^{-1},gh)
  \left(h^{-1}\cdot a_g(h\cdot v,h\cdot w)\right)
$$
in place of (\ref{eqn:conjugate}), and (\ref{eqn:jacobi}) is unchanged
(but now considered as a relation in $S(V)$). That is, these are
necessary and sufficient conditions for $A$ to be additively isomorphic
to $S(V)\ot \C G [t]$, leading to a generalization of Theorem \ref{thm:necsuff}.
}
\end{remark}

In the remainder of this section, we give several examples.

\begin{example}\label{ex:2}{\rm
(Symplectic groups.)
Let $V$ be a finite dimensional symplectic vector space over $\C$,
that is $\dim V$ is even and there is a nondegenerate skew-symmetric
form $\omega: V\times V\rightarrow \C$. Let $G$ be a finite subgroup
of the symplectic group {\em Sp}$(V)$ of all invertible linear
transformations preserving $\omega$. Let $\alpha: G\times G\rightarrow
\C^{\times}$ be the trivial two-cocycle $\alpha(g,h)=1$ for all
$g,h\in G$. For each element $g\in G$ such that $\codim V^g =2$
(the {\em symplectic reflections}), let $a_g=c_g\omega_g$, where $\omega_g$
is defined to be $\omega | _{(V^g)^{\perp}}$ on $(V^g)^{\perp}=\Ima (1-g)$
and $0$ on $V^g=\ker(1-g)$, and where $c_g\in \C$ are scalars such that
$c_{hgh^{-1}}=c_g$ for all $h\in G$. Let $a_1=c\omega$ for a scalar
$c$ and $a_g=0$ for all other $g\in G$. By 
\cite[Thm.\ 1.3]{etingof-ginzburg02}, the collection $\{a_g\mid g\in G\}$ 
determines a graded Hecke algebra, called a {\em symplectic reflection
algebra}. By Corollary \ref{cor:parameters},
the parameter space for the possible {\em twisted} graded Hecke
algebras cannot be larger than that for graded Hecke algebras in the
symplectic case (cf.\ \cite{chmutova}).
}\end{example}

The next two examples, by contrast, involve groups for which the parameter
space for twisted graded Hecke algebras is larger than that for graded
Hecke algebras.

\begin{example}\label{ex:1}
{\rm (Elementary abelian groups.)
Let $n\geq 3$, $\ell\geq 2$, $V=\C^n$ and $G\cong (\Z/\ell \Z)^{n-1}$
the multiplicative subgroup of $M_{n\times n}(\C)$ 
generated by the diagonal matrices
(where $q$ is a primitive $\ell$th root of 1):
\begin{eqnarray*}
  g_1&=&\diag(q,q^{-1},1,\ldots,1)\\
  g_2&=&\diag(1,q,q^{-1},1,\ldots,1)\\
     &\vdots & \\
  g_{n-1}&=&\diag(1,\ldots,1,q,q^{-1})
\end{eqnarray*}
Define a function $\alpha: G\times G\rightarrow \C^{\times}$ by
$$
  \alpha(g_1^{i_1}\cdots g_{n-1}^{i_{n-1}}, \ g_1^{j_1}\cdots
  g_{n-1}^{j_{n-1}}) = q^{-\sum_{1\leq k\leq n-2} i_kj_{k+1}}.
$$
The relation (\ref{eqn:alpha}) may be checked directly.
For each generator $g_i$ ($i=1,\ldots,n-1$), 
as well as for
$g_n=g_1^{-1}\cdots g_{n-1}^{-1}$, we have $\det(h^{\perp})=\alpha(h,g_i)
\alpha^{-1}(g_i,h)$ for all $h\in G$, as follows from a
straightforward computation. 
Thus we have a twisted graded Hecke algebra in which the forms
$a_{g_1},\ldots,a_{g_n}$ are all nonzero (and other $a_g=0$).
(Note that some of the relevant determinants
are not 1. There is no nontrivial (untwisted) graded Hecke algebra
in this case.) If $\ell\neq 2$, it may be checked that this gives the full
parameter space of twisted graded Hecke algebras corresponding to this
cocycle $\alpha$. If $\ell =2$, this determinant condition is in fact
satisfied for all $g$ with $\codim(V^g)=2$, and so there is a larger
parameter space of twisted graded Hecke algebras in this case.
This example is discussed in a different context in 
\cite[Example 4.1]{witherspoon}, 
which generalizes \cite[Example 4.7]{caldararu-giaquinto-witherspoon04}.
}\end{example}

\begin{example}\label{ex:3}{\rm 
(Symmetric groups.)
Let $G=S_n$, acting on $V=\C^n$ by the permutation representation.
Let $\alpha$ be the unique nontrivial two-cocycle (up to coboundary).
The Schur representation group for $S_n$ is
$$
  \Gamma_n=\langle t_1,\ldots,t_{n-1},\tau\mid\tau^2=1,
  t_i^2=1, \tau t_i=t_i\tau, (t_i t_{i+1})^3=1,
  (t_it_j)^2=\tau (i\leq j-2)\rangle.
$$
(See for example \cite[p.\ 179]{karpilovsky85}.)
A two-cocycle $S_n\times S_n\rightarrow \C^{\times}$
is determined up to coboundary by an irreducible representation of the Schur 
representation group $\Gamma_n$ for which the central
element $\tau$ necessarily acts as multiplication by a scalar. As $\tau^2=1$,
this scalar must be $-1$ in case of the nontrivial cocycle
$\alpha$. Let $g=(12)(34)$. The presentation of $\Gamma_n$
may be used to check that 
$$
  \det(h^{\perp})=\alpha(h,g)\alpha^{-1}(g,h)
$$
for all $h\in C(g)$.
For example, letting $h=(12)$, we have $\det(h^{\perp})=-1$.
Choose a section $T:S_n\rightarrow \Gamma_n$ of the projection from 
$\Gamma_n$ to $S_n$ that sends $\tau$ to $1$ and $t_i$ to the 
transposition $(i, i+1)$ as follows: Let $T(12)=t_1$, $T(34)=t_3$,
$T((12)(34))=t_1t_3$, and choose other images arbitrarily. 
Then $\alpha(h,g)\alpha^{-1}(g,h)$ is the
scalar by which the following element acts on an irreducible
representation of $\Gamma_n$:
$$
\hspace{-.3in}T(12)T((12)(34))T((12)(12)(34))^{-1} T((12)(34)(12))
  T(12)^{-1}T((12)(34))^{-1}
$$

\vspace{-.28in}

\begin{eqnarray*}
   \hspace{1in}&=& T(12)T(12)T(34)T(34)^{-1} T(34)T(12)^{-1}T((12)(34))^{-1}\\
  \hspace{1in}&=& T(34)T(12)T(34)T(12)\\
  \hspace{1in}&=& \tau,
\end{eqnarray*}
thus the scalar is $-1$ as desired. 
Other elements $h$ of $C(g)$ may be checked similarly; as
$h\mapsto \det(h^{\perp})$ is a group homomorphism, it suffices to
check the condition (\ref{eqn:condition}) on generators of $C(g)$.
Therefore there
is a twisted graded Hecke algebra corresponding to $\alpha$
with $a_g\neq 0$. (Compare with trivial $\alpha$, where $a_g$
is necessarily 0 for this choice of $g$ \cite[Table 1]{ram-shepler03}.)
With the choice $g=(123)$, we find as in the case of the trivial
cocycle $\alpha$, that $1=\det(h^{\perp})=\alpha(h,g)\alpha^{-1}(g,h)$
for all $h\in C(g)$. Therefore the parameter space of twisted graded
Hecke algebras is larger than that of graded Hecke algebras,
and involves the conjugacy classes of 
all $g\in G$ for which $\codim V^g =2$.
}\end{example}

\section{Deformations of crossed products}

In this section, we prove that the twisted graded Hecke algebras
of the previous section are precisely the deformations of
$S(V)\#_{\alpha}G$ of a particular type. 
We use this connection to deformations
to put our results and examples in a larger context.

Let $t$ be an indeterminate.
Given any associative algebra $R$ over $\C$ (for example
$R=S(V)\#_{\alpha}G$), a {\em deformation of $R$ over $\C[t]$} is an
associative $\C[t]$-algebra with underlying vector space 
$R[t]=\C[t]\ot R$ and multiplication of the form
$$
  r*s=rs+\mu_1(r,s)t +\mu_2(r,s)t^2+\cdots +\mu_p(r,s)t^p
$$
for all $r,s\in R$, where $rs$ denotes the product of $r$ and $s$
in $R$, the $\mu_i:R\times R\rightarrow R$ are $\C$-bilinear
maps extended to be $\C[t]$-bilinear, and $p$ depends on $r,s$. 
Associativity of $*$ implies
that $\mu_1$ is a {\em Hochschild two-cocycle}, that is
\begin{equation}\label{eqn:hochschild}
  \mu_1(w,r)s+\mu_1(wr,s)=\mu_1(w,rs)+w\mu_1(r,s)
\end{equation}
for all $w,r,s\in R$, as well as further conditions on the
$\mu_i$, $i\geq 1$. 
The cocycle $\mu_1$ is called the {\em infinitesimal} of the
deformation. 
(See \cite{gerstenhaber64} or \cite{giaquinto-zhang98}
for the details from algebraic deformation theory.)

Recall that $S(V)\#_{\alpha}G$ is a graded algebra where we assign
degree 1 to elements of $V$ and degree 0 to elements of $G$. 

\begin{thm}\label{thm:equivalence}
Up to isomorphism,
the twisted graded Hecke algebras are precisely the
deformations of $S(V)\#_{\alpha}G$ over $\C[t]$ 
for which $\deg \mu_i = -2i$ ($i\geq 1$).
\end{thm}

\begin{proof}
Let $A$ be a twisted graded Hecke algebra, as defined in 
Section 2. 
Let $v_1,\ldots,
v_n$ be a basis of $V$ and choose the section of the projection 
from $T(V)$ to $S(V)$ in which a word in $T(V)$ is written in the
order $v_1^{i_1}\cdots v_n^{i_n}$.
Express all elements of $A$ in terms of this
section, writing group elements on the right. Such expressions in $A$
exist and are unique due to the additive isomorphism $A\cong
S(V)\#_{\alpha}G[t]$. Now let $r=v_1^{i_1}\cdots v_n^{i_n}
\overline{g}$ and $s=v_1^{j_1}\cdots v_n^{j_n}\overline{h}$ be elements of $A$.
Denoting the product in $A$ by $*$, we have
$$
  r*s=\alpha(g,h)v_1^{i_1}\cdots v_n^{i_n}*(g\cdot v_1^{j_1}
  \cdots v_n^{j_n})\overline{gh}.
$$
The factor $v_1^{i_1}\cdots v_n^{i_n}(g\cdot v_1^{j_1}\cdots
v_n^{j_n})$ may now be rearranged using 
(\ref{eqn:bracket}) repeatedly until it is in the standard form
for elements of $A$ discussed above.
The result will be of the form
$$
  r*s=rs+\mu_1(r,s)t+\mu_2(r,s)t^2+\cdots +\mu_p(r,s)t^p,
$$
where $rs$ is identified with the product of $r$ and $s$ in $S(V)\#_{\alpha}G$,
via the additive isomorphism $A\cong S(V)\#_{\alpha}G[t]$,
and $\mu_1(r,s),\ldots,\mu_p(r,s)$ are elements of $A$ also identified with
elements of $S(V)\#_{\alpha}G$.
This is a finite process as each time (\ref{eqn:bracket}) is applied,
the degree drops.
As the multiplication $*$ in the twisted graded Hecke algebra
$A$ is bilinear and associative, the maps $\mu_i$ are bilinear
and thus $A$ is a deformation of $S(V)\#_{\alpha}G$ over $\C[t]$.
(Alternatively, the existence of the obvious algebra isomorphism 
$A/tA\cong S(V)\#_{\alpha}G$ implies
that $A$ is a deformation of $S(V)\#_{\alpha}G$ over $\C[t]$.)
The conditions on the $\mu_i$ stated in the theorem are consequences of the
relations (\ref{eqn:bracket}), by induction on the degree 
$\sum_{k=1}^n (i_k+j_k)$ of a 
product $v_1^{i_1}\cdots v_n^{i_n}*v_1^{j_1}\cdots v_n^{j_n}$.

Conversely, suppose that $A$ is a deformation of
$S(V)\#_{\alpha}G$ over $\C[t]$ satisfying the given conditions.
By definition, $A\cong S(V)\#_{\alpha}G[t]$ as a vector space
over $\C[t]$. Define a $\C[t]$-linear map 
$\phi: T(V)\#_{\alpha}G[t]\rightarrow A$ by
$$
  \phi(v_{i_1}\cdots v_{i_m}\overline{g}) = v_{i_1} * \cdots
    * v_{i_m} *\overline{g}
$$
for all words $v_{i_1}\cdots v_{i_m}$ in $T(V)$ and $g\in G$.
It may be checked that 
$\phi$ is an algebra homomorphism since $T(V)$ is free on
$v_1,\ldots, v_n$, and since $\mu_i(\C^{\alpha}G,\C^{\alpha}G)=
\mu_i(\C^{\alpha}G,V)=\mu_i(V,\C^{\alpha}G)=0$ for all $i\geq 1$ by
the degree condition on $\mu_i$.
We will show, by induction on degree, that $\phi$ is surjective: 
First note that 
$\phi(\overline{g})=\overline{g}$ and $\phi(v\overline{g})=v
\overline{g}$ for all $v\in V$, $g\in G$.
Now we would like to show that an arbitrary basis monomial 
$v_{i_1}\cdots v_{i_m}\overline{g}$ ($i_1\leq\cdots\leq i_m$)
of $A$ is in $\Ima(\phi)$.
Assume $v_{i_2}\cdots v_{i_m}\overline{g}=\phi(X)$ for some element
$X$ of $T(V)\#_{\alpha}G[t]$. Then 
\begin{eqnarray*}
\phi(v_{i_1}X)&=&\phi(v_{i_1})*\phi(X) \\
 &=& v_{i_1} * v_{i_2}\cdots v_{i_m}\overline{g}\\
 &=& v_{i_1}\cdots v_{i_m}\overline{g} + \mu_1(v_{i_1},v_{i_2}\cdots v_{i_m}
    \overline{g})t + \mu_2(v_{i_1},v_{i_2}\cdots v_{i_m}
    \overline{g})t^2 + \cdots
\end{eqnarray*}
By induction on $m$, as $\mu_i$ is a map of degree $-2i$,
each $\mu_j(v_{i_1},v_{i_2}\cdots v_{i_m}\overline{g})$ is in $\Ima (\phi)$.
This implies $v_{i_1}\cdots v_{i_m}\overline{g}\in \Ima (\phi)$,
and thus $\phi$ is surjective.

It remains to determine the kernel of $\phi$.
Letting $v,w\in V$, we find that
\begin{eqnarray*}
  \phi(vw) & = & v*w =vw+\mu_1(v,w)t\\
   \phi(wv) &=& w*v=wv+\mu_1(w,v)t,
\end{eqnarray*}
since $\deg \mu_i = -2i$.
As $vw=wv$ in $S(V)$, we may subtract to obtain $\phi(vw-wv)=
(\mu_1(v,w)-\mu_1(w,v))t$. Since $\deg \mu_1=-2$ and $\phi
(\overline{g})=\overline{g}$ for all $g\in G$, this implies that
\begin{equation}\label{eqn:kernel}
  vw-wv-(\mu_1(v,w)-\mu_1(w,v))t
\end{equation}
is in the kernel of $\phi$
for all $v,w\in V$. It also follows, as $\deg\mu_1 = -2$, that 
$$
  \mu_1(v,w)-\mu_1(w,v)=\sum_{g\in G} a_g(v,w)\overline{g}
$$
for some functions $a_g: V\times V\rightarrow \C$. By
definition, $a_g$ is bilinear and skew-symmetric for each $g\in G$.
Let $I[t]$ be the ideal of $T(V)\#_{\alpha}G[t]$ generated by
all such expressions (\ref{eqn:kernel}), so that $I[t]\subset\Ker \phi$.
We will use a dimension count to show that $I[t]=\Ker \phi$:
By the arguments of the previous section, $T(V)\#_{\alpha}G[t]/I[t]$
is a quotient of $S(V)\#_{\alpha}G[t]$, and so has dimension
in each degree no greater than that of $S(V)\#_{\alpha}G[t]$.
Since $\phi$ induces a map from $T(V)\#_{\alpha}G[t]/I[t]$ onto
$S(V)\#_{\alpha}G[t]$, this forces $I[t]=\Ker\phi$ and thus $A$
is a twisted graded Hecke algebra.
\end{proof}

The proof of Theorem \ref{thm:equivalence} shows that the
skew-symmetric forms $a_g$ appearing in the relations
(\ref{eqn:bracket}) of a twisted graded Hecke algebra arise
as coefficients in the skew-symmetrization of a Hochschild 
two-cocycle $\mu_1$ on $S(V)\#_{\alpha}G$. Hochschild
cohomology thus provides an alternate approach to the study
of twisted graded Hecke algebras. An advantage of this
approach is that it puts the computation of the possible
skew-symmetric forms $a_g$ into a larger context. A 
disadvantage is that, given a Hochschild two-cocycle $\mu_1$
on an algebra $R$, there is no general method for determining
whether it lifts to a deformation of $R$, nor for
finding the $\mu_i$ ($i\geq 2$) in case it does lift. 
In case $\deg(\mu_1)=-1$, Remark \ref{deg1} points to existence
of deformations of $S(V)\#_{\alpha}G$ defined by quotients
of $T(V)\#_{\alpha}G[t]$ analogous to twisted graded Hecke algebras.
There are a few special cases where deformations of $S(V)\#_{\alpha}G$
are known whose infinitesimals $\mu_1$ have arbitrarily high degree
\cite{caldararu-giaquinto-witherspoon04,witherspoon}.
These examples were discovered after examining the relevant Hochschild
cohomology:
$$
  \HH^{\DOT} (S(V)\#_{\alpha}G)\cong \left(\bigoplus_{g\in G}
  \Wedge^{\DOT -\codim V^g}((V^g)^*) \ot \det(((V^g)^{\perp})^*)\ot
   S(V^g)\overline{g}\right) ^G
$$
additively, where the superscript $G$ denotes elements invariant under the
action induced by conjugation in $\C^{\alpha}G$ and the action
of $G$ on $V$, and $\det$ denotes the top exterior power.
This follows from \cite{ginzburg-kaledin04}, see (6.2) and the formula
before (6.4) where $S(V)$ is replaced by $\C[M]$, or more directly
from the techniques of \cite{AFLS}, where $S(V)$ is replaced by
a Weyl algebra.

\quad

We conclude by giving a broader context for the
isomorphism found by Ram and Shepler between Lusztig's and Drinfeld's
definitions of (untwisted) graded Hecke algebras for real reflection groups
\cite{drinfeld86,lusztig89,ram-shepler03}.
Let $W\subset GL(V)$ 
be a finite real reflection group. Specifically, $W$ is generated by
simple reflections $s_1,\ldots,s_n$ corresponding to simple roots
$\alpha_1,\ldots,\alpha_n\in V$ for a root system $R$ in $V$.
If $\alpha\in R$ is any root, we write $s_{\alpha}$ for
the reflection corresponding to $\alpha$, so that
$$
  s_{\alpha}\cdot v = v-\langle v,\alpha\check{\hspace{.2cm}}\rangle\alpha
$$   
where $\alpha\check{\hspace{.2cm}}= 2\alpha / \langle
  \alpha , \alpha\rangle$, and $\langle \ , \ \rangle$ is
the inner product.
For each $\alpha\in R$, choose $k_{\alpha}\in\C$ in such a way that
$k_{g\cdot\alpha}=k_{\alpha}$ for all $g\in W$, so that
the number of independent parameters $k_{\alpha}$ is simply the number
of different lengths of roots. 
Lusztig's graded version of an affine Hecke algebra, 
as defined over $\C[t]$,
is the quotient $A'$ of $T(V)\# W[t]$ by the
ideal generated by all 
$$
  [v,w] \ \ \mbox{ and } \ \ \overline{s}_iv-(s_i\cdot v)\overline{s}
  _i + k_{\alpha_i}
  \langle v,\alpha_i\check{\hspace{.2cm}}\rangle t
$$
where $v,w\in V$ and $i=1,\ldots, n$.

Ram and Shepler defined the following graded Hecke algebra and showed
that it is isomorphic to Lusztig's algebra $A'$ given above (and their
isomorphism extends to one over $\C[t]$): For each $g\in W$, let
$$
  a_g(v,w)=\frac{1}{4}\sum_{\stackrel{\alpha,\beta>0}
  {g=s_{\alpha}s_{\beta}}} k_{\alpha}k_{\beta} \left(
  \langle v,\beta\check{\hspace{.2cm}}\rangle\langle w,
  \alpha\check{\hspace{.2cm}}\rangle -\langle v,
  \alpha\check{\hspace{.2cm}}\rangle \langle w,
  \beta\check{\hspace{.2cm}}\rangle\right).
$$
In order to make the connection to Lusztig's version, replace
$t$ by $t^2$ in (\ref{eqn:bracket}) and take the quotient $A$ of $T(V)\# W[t]$
by the ideal generated by all 
$$
  [v,w] - \sum_{g\in W}a_g(v,w)t^2 \overline{g}.
$$
The resulting quotient algebra $A$ is in fact a graded Hecke algebra
defined over $\C[t^2]$, however
the isomorphism with Lusztig's algebra $A'$ is not defined
over $\C[t^2]$. This isomorphism
$\Phi_t : A'\rightarrow A$ is given by
$$
  \Phi_t (v) = v- t\langle v,h\rangle \ \ \mbox{ and } \ \
  \Phi_t(\overline{g})= \overline{g}
$$
for all $v\in V, g\in W$, where $h=\displaystyle{\frac{1}{2}}
\sum_{\alpha >0}k_{\alpha}\alpha{\check{\hspace{.2cm}}}\overline{s}_{\alpha}$,
and $\langle v,h\rangle = \displaystyle{\frac{1}{2}\sum_{\alpha >0}
k_{\alpha}\langle v,\alpha\check{\hspace{.2cm}}\rangle 
\overline{s}_{\alpha}}$. 
This may be extended uniquely to an algebra
isomorphism by the calculations of Ram and Shepler.
The reason for
replacing $t$ by $t^2$ in the definition of $A$ is that it is
necessary in order to extend their isomorphism to one defined over
$\C[t]$. The only calculation that changes significantly in this context
occurs in the proof that
$$
  [\Phi_t(v),\Phi_t(w)] = [v,w] + t^2[\langle v,h\rangle , 
  \langle w,h\rangle ] -t[v,\langle w,h\rangle ] +t
  [w,\langle v,h\rangle ] = 0.
$$
This is true as
$[v,\langle w,h\rangle ] = [w,\langle v,h\rangle ]$ and
$[\langle v,h\rangle , \langle w,h\rangle ] = -\sum_{g\in W}a_g(v,w)
\overline{g}$ by computations in the proof of 
\cite[Thm.\ 3.5]{ram-shepler03}.

In particular, the isomorphism $\Phi_t$ is an equivalence of 
deformations over $\C[t]$.
This implies that the Hochschild two-cocycle corresponding to the
coefficients of $t$ in the deformation
$A'$ of $S(V)\# W$ is a coboundary, however there is a Hochschild
two-cocycle arising from the coefficients of $t^2$ in 
products in $A'$ that
will be cohomologous to that for $A$. (Again see 
\cite{gerstenhaber64} or \cite{giaquinto-zhang98} for details from the general 
theory of algebraic deformations.)

\end{document}